\documentclass[a4paper,12pt]{article}

\makeatletter
\@addtoreset{footnote}{page}
\makeatother

\usepackage{amsthm}
\usepackage{amsmath,amssymb,latexsym,amsfonts,mathrsfs}
\usepackage{bm}

\newcommand{\eps}{\ensuremath{\varepsilon}}
 
\newcommand{\dist}{\stackrel{{\rm d}}{=}}
\newcommand{\tend}[2]{\mathrel{\mathop{\longrightarrow}\limits^{#1}_{#2}}}
\renewcommand{\hat}{\widehat}
\renewcommand{\tilde}{\widetilde}
\renewcommand{\bar}{\overline}

\newcommand{\rbra}[1]{\!\left( #1 \right)} 
\newcommand{\cbra}[1]{\!\left\{ #1 \right\}} 
\newcommand{\sbra}[1]{\!\left[ #1 \right]} 

\newcommand{\bE}{\ensuremath{\mathbb{E}}}
\newcommand{\bN}{\ensuremath{\mathbb{N}}}
\newcommand{\bP}{\ensuremath{\mathbb{P}}}
\newcommand{\bR}{\ensuremath{\mathbb{R}}}
\newcommand{\cD}{\ensuremath{\mathcal{D}}}
\newcommand{\cF}{\ensuremath{\mathcal{F}}}
\newcommand{\cM}{\ensuremath{\mathcal{M}}}
\newcommand{\cN}{\ensuremath{\mathcal{N}}}
\newcommand{\cS}{\ensuremath{\mathcal{S}}}
\newcommand{\vA}{\ensuremath{\bm{A}}}
\newcommand{\vB}{\ensuremath{\bm{B}}}
\newcommand{\vI}{\ensuremath{\bm{I}}}
\newcommand{\vP}{\ensuremath{\bm{P}}}

\newcommand{\vR}{\ensuremath{\bm{R}}}
\newcommand{\vX}{\ensuremath{\bm{X}}}
\newcommand{\vY}{\ensuremath{\bm{Y}}}
\newcommand{\vZ}{\ensuremath{\bm{Z}}}
\newcommand{\va}{\ensuremath{\bm{a}}}
\newcommand{\vb}{\ensuremath{\bm{b}}}

\newcommand{\vq}{\ensuremath{\bm{q}}}
\newcommand{\vx}{\ensuremath{\bm{x}}}

\theoremstyle{plain}
\newtheorem{Thm}{Theorem}[section]
\newtheorem{Lem}[Thm]{Lemma}

\theoremstyle{definition}

\newcommand{\Proof}[2][Proof]{\begin{proof}[{#1}] #2 \end{proof}}

\numberwithin{equation}{section}

\makeatletter
\renewcommand\section{\@startsection {section}{1}{\z@}%
                                   {-3.5ex \@plus -1ex \@minus -.2ex}%
                                   {2.3ex \@plus.2ex}%
                                   {\normalfont\large\bf}}
\makeatother

\makeatletter
\renewcommand\subsection{\@startsection {subsection}{1}{\z@}%
                                   {-3.5ex \@plus -1ex \@minus -.2ex}%
                                   {2.3ex \@plus.2ex}%
                                   {\normalfont\normalsize\bf}}
\makeatother

\newcommand{\vpsi}{\ensuremath{\bm{\psi}}}
\newcommand{\vdelta}{\ensuremath{\bm{\delta}}}
\newcommand{\vdI}{\ensuremath{\bm{\delta I}}}
\newcommand{\vnu}{\ensuremath{\bm{\nu}}}
\newcommand{\vzero}{\ensuremath{\bm{0}}}

\begin{document}

\begin{center}
{\Large \bf 
Asymptotic normality in multi-dimension of nonparametric estimators 
for discrete-time semi-Markov chains
}
\end{center}
\footnotetext{
This research was supported by RIMS and by ISM.}
\begin{center}
Hiroki Ogata\footnote{
Graduate School of Science, Kyoto University.}, 
\quad 
Luis Iv\'an Hern\'andez Ru\'{\i}z\footnotemark[1]\footnote{
The research of this author was supported 
by JSPS Open Partnership Joint Research Projects grant no. JPJSBP120209921. 
} 
\quad 
and 
\quad 
Kouji Yano\footnotemark[2]\footnote{
Graduate School of Science, Osaka University.}\footnote{
The research of this author was supported by 
JSPS KAKENHI grant no.'s JP19H01791, JP19K21834 and 21H01002.} 
\end{center}

\begin{abstract}
The asymptotic normality in multi-dimension of the nonparametric estimator 
of the transition probabilities of a Markov renewal chain is proved, 
and is applied to that of other nonparametric estimators involved 
with the associated semi-Markov chain. 
\end{abstract}

\section{Introduction}

In the case of continuous time, 
asymptotic normality of the nonparametric estimator 
for the transition probabilities of a continuous-time Markov renewal chain 
has been studied by Moore--Pyke \cite{MR240934} in one dimension, 
and was generalized by Gill \cite{MR576899} to multi-dimension. 
Ouhbi--Limnios \cite{MR1918880,MR1946645} studied 
the asymptotic normality of other nonparametric estimators involved 
with the associated semi-Markov chain. 
Several statistical properties of those estimators 
were studied in \cite{MR513932,MR574546,MR933988,MR1417674,MR2139370,MR2163098,MR2245576,MR4423721}, 
and some applications were discussed in \cite{Heutte2002}. 

In the case of discrete time, 
the asymptotic normality in one-dimension 
of various nonparametric estimators has been investigated by 
Barbu--Limnios \cite{BL}. 
Our aim is to establish the asymptotic normality in multi-dimension 
of the nonparametric estimator of the transition probabilities 
of the underlying Markov renewal chain, and to provide a systematic and greatly simplified method 
of deriving those results of the asymptotic normality by Barbu--Limnios \cite{BL}. 

Let $ E $ be a finite state space 
and $ \bN = \{ 0,1,2,\ldots \} $ denote the set of nonnegative integers. 
We write $ \tilde{E} := E \times E \times \bN $. 
Let $ \vq = (\vq(k))_{k \in \bN} = (q_{ij}(k))_{(i,j,k) \in \tilde{E}} $ 
be a matrix sequence with nonnegative entries 
such that $ q_{ij}(0) = 0 $ for $ i,j \in E $ and 
\begin{align}
\sum_{j \in E} \sum_{k \in \bN} q_{ij}(k) = 1 
\quad \text{for all $ i \in E $} . 
\label{}
\end{align}
A \emph{Markov renewal chain} is a time-homogeneous Markov chain 
$ (J,S) = (J_n,S_n)_{n=0}^{\infty } $ with values in $ E \times \bN $ 
defined on probability spaces $ (\Omega,\cF,\bP_{i_0}) $ for $ i_0 \in E $ such that 
the process $ (J,X) = (J_n,X_n)_{n=0}^{\infty } $ 
defined by $ J_0=i_0 $, $ X_0 := S_0 := 0 $ 
and $ X_n := S_n-S_{n-1} $ for $ n \ge 1 $ is a Markov chain 
whose transition probabilities are given as 
\begin{align}
\bP_{i_0}(J_{n+1} = j , \ X_{n+1}= k \mid J_n = i , \ X_n) = q_{ij}(k) 
\quad \text{for all $ i,j \in E $ and $ n,k \in \bN $} . 
\label{}
\end{align}
(The quantity $ q_{ij}(k) $ is often called the \emph{semi-Markov kernel}, 
although it describes the transition of the Markov renewal chain, 
not of the associated semi-Markov chain.) 
Note that the \emph{embedded Markov chain} $ J = (J_n)_{n=0}^{\infty } $ 
is an $ E $-valued Markov chain whose transition probabilities are given as 
\begin{align}
\bP_{i_0}(J_{n+1} = j \mid J_n = i) = q^J_{ij} := \sum_{k \in \bN} q_{ij}(k) 
\quad \text{for all $ i,j \in E $ and $ n \in \bN $} . 
\label{}
\end{align}

For $ j \in E $, we write $ \tau_j = \inf \{ n \ge 1 : J_n = j \} $. 
We assume the following conditions: 
\begin{description}
\item[(A1)] 
The embedded Markov chain $ J $ is irreducible, i.e., 
$ \bP_i(\tau_j < \infty ) > 0 $ for all $ i,j \in E $. 
\item[(A2)] 
The Markov renewal chain $ (J,S) $ is aperiodic, i.e., 
the greatest common divisor of $ \{ k \in \bN : \bP_i(S_{\tau_i} = k) > 0 \} $ is one 
for all $ i \in E $. 
\item[(A3)] 
The Markov renewal chain $ (J,S) $ is positive recurrent, i.e., 
$ \mu_{ii} := \bE_i S_{\tau_i} < \infty $ for all $ i \in E $. 
\end{description}
We write {\bf (A)} for all of the above assumptions being satisfied. 
In this case, it holds (see, e.g., \cite[Theorem 1.7.7]{N}) 
that the embedded Markov chain $ J $ 
has a unique invariant probability $ \vnu = (\nu_i)_{i \in E} $, 
and it also holds (see Proposition 3.6 of \cite{BL}) that 
\begin{align}
\bar{m} := \sum_{i \in E} \nu_i \bE_i S_1 = \mu_{jj} \nu_j 
\quad \text{for all $ j \in E $}. 
\label{eq: barm}
\end{align}

For $ M \in \bN $ and for $ i \in E $, we set 
\begin{align}
N(M) := \max \{ n \in \bN : S_n \le M \} 
\quad \text{and} \quad 
N_i(M) := \sum_{n=1}^{N(M)} 1_{\{ J_{n-1}=i \}} 
\label{}
\end{align}
so that $ N(M) = \sum_{i \in E} N_i(M) $. 
For the unknown quantity $ \vq $, 
we consider the nonparametric estimator 
\begin{align}
\hat{\vq}(M) = (\hat{\vq}(k,M))_{k \in \bN} 
= (\hat{q}_{ij}(k,M))_{(i,j,k) \in \tilde{E}} 
\label{eq: hatvq}
\end{align}
which is defined as 
\begin{align}
\hat{q}_{ij}(k,M) 
:=& \frac{1}{N_i(M)}\sum_{n=1}^{N(M)} 1_{\{J_{n-1}=i, \ J_{n}=j, \ X_{n}=k\}} 
\label{}
\end{align}
so that $ \sum_{j \in E} \sum_{k \in \bN} \hat{q}_{ij}(k,M) = 1 $ for all $ i \in E $. 
The strong consistency of $ \hat{\vq}(M) $ is known as follows. 

\begin{Thm}[Barbu--Limnios {\cite[Theorem 4.1]{BL}}]
Suppose {\bf (A)} is satisfied. Then 
\begin{align}
\hat{\vq}(M) \tend{\rm a.s.}{M \to \infty } \vq , 
\label{}
\end{align}
or equivalently, 
\begin{align}
\hat{q}_{ij}(k,M) 
\tend{\rm a.s.}{M \to \infty } q_{ij}(k) 
\quad \text{for all $ (i,j,k) \in \tilde{E} $}, 
\label{}
\end{align}
under $ \bP_{i_0} $ for all initial value $ i_0 \in E $. 
\end{Thm}

Let $ \cN(0,v) $ stand for the centered normal distribution whose variance is $ v $. 
The asymptotic normality of $ \hat{\vq}(M) $ in one-dimension 
has been obtained as follows. 

\begin{Thm}[Barbu--Limnios {\cite[Theorem 4.2]{BL}}] \label{qone}
Suppose {\bf (A)} is satisfied. 
Then, for any fixed $ (i,j,k) \in \tilde{E} $, it holds that 
\begin{align}
\sqrt{M} \rbra{ \hat{q}_{ij}(k,M) -q_{ij}(k) \Big. } 
\tend{\rm d}{M \to \infty } 
\cN(0,v^{\vq}(i,j,k)) 
\label{}
\end{align}
under $ \bP_{i_0} $ for all initial value $ i_0 \in E $, 
where the asymptotic variance $ v^{\vq}(i,j,k) $ is given as 
\begin{align}
v^{\vq}(i,j,k) := \mu_{ii} q_{ij}(k) (1-q_{ij}(k)) . 
\label{}
\end{align}
\end{Thm}

Let us discuss convergence in multi-dimension. 
For a nonnegative definite matrix $ V = (V_{e e'})_{e,e' \in \tilde{E}} $, 
we write $ \cN(\vzero,V) $ for the distribution of 
a centered Gaussian system $ Y=(Y_{e})_{e \in \tilde{E}} $ 
whose covariances are given as $ \bE Y_{e} Y_{e'} = V_{e e'} $, 
and we call $ V $ the \emph{variance matrix} of $ \cN(\vzero,V) $. 
We now state our main theorem about 
the asymptotic normality of $ \hat{\vq}(M) $ in multi-dimension. 

\begin{Thm} \label{thm: q}
Suppose {\bf (A)} is satisfied. 
Then, on the product space 
\begin{align}
\bR^{\tilde{E}} = \cbra{ \vx = (x_e)_{e \in \tilde{E}} 
= (x_{ij}(k))_{(i,j,k) \in \tilde{E}} : 
\text{$ x_e \in \bR $ for $ e \in \tilde{E} $} } 
\label{}
\end{align}
equipped with the product topology, 
it holds that 
\begin{align}
\sqrt{M} \rbra{ \hat{\vq}(M) - \vq \Big. } 
\tend{\rm d}{M \to \infty } 
\cN(\vzero,V^{\vq}) 
\label{eq: q}
\end{align}
under $ \bP_{i_0} $ for all initial value $ i_0 \in E $, 
where the asymptotic variance matrix $ V^{\vq} = (V^{\vq}_{e e'})_{e,e' \in \tilde{E}} $ 
is given as 
\begin{align}
V^{\vq}_{(i,j,k) (i',j',k')} 
= \mu_{ii} 1_{\{ i=i' \}} q_{ij}(k) \rbra{ 1_{\{ j=j', \ k=k' \}} - q_{i'j'}(k') } . 
\label{}
\end{align}
\end{Thm}

Noting that $ V^{\vq}_{(i,j,k)(i,j,k)} = v^{\vq}(i,j,k) $, 
we can derive the result of Theorem \ref{qone} 
directly from our Theorem \ref{thm: q}.

Let us consider 
the \emph{semi-Markov chain} $ Z = (Z_k)_{k \in \bN} $ associated to $ (J,S) $, 
which is defined as the time-changed process $ Z_k := J_{N(k)} $ for $ k \in \bN $ with 
$ N(k) := \max \{ n \in \bN : S_n \le k \} $. 
Note that $ Z_0 = J_0 $. 
We set 
\begin{align}
P_{ij}(k) := \bP_i(Z_k = j) , 
\label{eq: P}
\end{align}
and we call $ \vP = (\vP(k))_{k \in \bN} = (P_{ij}(k))_{(i,j,k) \in \tilde{E}} $ 
the \emph{distribution matrix sequence} of the semi-Markov chain $ Z $. 

Let $ E = U \cup D $ be a nontrivial partition of the finite state space $ E $ 
and write 
\begin{align}
T_D = \inf \{ k \in \bN : Z_k \in D \} . 
\label{}
\end{align}
We define 
\begin{align}
R_i(k) = \bP_i(T_D>k) 
\quad \text{for $ i \in U $},  
\label{}
\end{align}
and call $ \vR = (\vR(k))_{k \in \bN} = (R_i(k))_{(i,k) \in U \times \bN} $ 
the \emph{reliability vector sequence}. 

As applications of our Theorem \ref{thm: q}, 
we shall generalize the asymptotic normality results in one-dimension to multi-dimension 
of certain estimators $ \hat{\vP}(M) $ and $ \hat{\vR}(M) $ 
of the distribution matrix sequence $ \vP $ and the reliability vector sequence $ \vR $, 
respectively.

This paper is organized as follows. 
In Section \ref{sec: MRC}, we prove our Theorem \ref{thm: q} 
about the asymptotic normality of $ \hat{\vq}(M) $. 
In Section \ref{sec: MCI}, we discuss matrix convolution inverse. 
Sections \ref{sec: DMS} and \ref{sec: RVS} are devoted to 
the asymptotic normality of the estimators $ \hat{\vP}(M) $ and $ \hat{\vR}(M) $, respectively. 
In Appendix, we deduce the asymptotic variances for the asymptotic normality in one dimension 
from our results in multi-dimension.

\section{Estimation of the transition probabilities} \label{sec: MRC} 

We write $ \va \cdot \vb $ for the inner product in $ \bR^{\tilde{E}} $: 
for $ \va = (a_e)_{e \in \tilde{E}} $ and $ \vb = (b_e)_{e \in \tilde{E}} $, set 
\begin{align}
\va \cdot \vb := \sum_{e \in \tilde{E}} a_e b_e . 
\label{}
\end{align}
We denote 
\begin{align}
\bR_0^{\tilde{E}} := \cbra{ \vx = (x_e)_{e \in \tilde{E}} \in \bR^{\tilde{E}} 
: \text{$ x_e = 0 $ for all but finite $ e \in \tilde{E} $} } 
\label{}
\end{align}
It is well-known (see, e.g., \cite[Theorem 4.29]{K}) that, 
for $ \bR^{\tilde{E}} $-valued processes, 
convergence in distribution on $ \bR^{\tilde{E}} $, 
namely 
\begin{align}
\vX(M) = (X_e(M))_{e \in \tilde{E}} 
\tend{\rm d}{M \to \infty } \vX = (X_e)_{e \in \tilde{E}} , 
\label{}
\end{align}
is equivalent to convergence of the finite-dimensional distributions, 
namely 
\begin{align}
(X_e(M))_{e \in I} 
\tend{\rm d}{M \to \infty } (X_e)_{e \in I} 
\label{}
\end{align}
for all finite subset $ I $ of $ \tilde{E} $. 
By the Cram\'er--Wold device, this convergence is reduced to the convergence 
\begin{align}
\va \cdot \vX(M) \tend{\rm d}{M \to \infty } \va \cdot \vX 
\label{}
\end{align}
for all constant $ \va = (a_e)_{e \in \tilde{E}} \in \bR_0^{\tilde{E}} $.

Let us proceed to the proof of our main theorem. 

\Proof[Proof of Theorem \ref{thm: q}]{
In order to obtain the desired convergence 
\begin{align}
\tilde{\vq}(M) := \sqrt{M} (\hat{\vq}(M) - \vq) \tend{\rm d}{M \to \infty } \cN(\vzero,V^{\vq}) , 
\label{}
\end{align}
it suffices to show that the convergence 
\begin{align}
\va \cdot \tilde{\vq}(M) 
\tend{\rm d}{M \to \infty } 
\cN(0,\va \cdot V^{\vq} \va) 
\label{eq: vatildevq}
\end{align}
holds for all constant $ \va = (a_{e})_{e \in \tilde{E}} \in \bR_0^{\tilde{E}} $, 
where we define 
\begin{align}
\va \cdot V^{\vq} \va 
:= \sum_{e,e' \in \tilde{E}} a_{e} V^{\vq}_{e e'} a_{e'} . 
\label{}
\end{align}
We may represent 
\begin{align}
\va \cdot \tilde{\vq}(M) 
= \vA(M) \cdot \vY(M) 
\label{}
\end{align}
if we set $ \vA(M) = (A_e(M))_{e \in \tilde{E}} $ 
and $ \vY(M) = (Y_e(M))_{e \in \tilde{E}} $ as 
\begin{align}
A_{(i,j,k)}(M) :=& a_{(i,j,k)} \frac{\sqrt{M N(M)}}{N_i(M)} , 
\label{} \\
Y_{(i,j,k)}(M) 
:=& \frac{1}{\sqrt{N(M)}} 
\sum_{n=1}^{N(M)} 1_{\{J_{n-1}=i\}} \rbra{ 1_{\{ J_{n}=j, \ X_{n}=k \}} - q_{ij}(k) } . 
\label{}
\end{align}
By \eqref{eq: barm} and by Propositions 3.7 and 3.8 of \cite{BL}, we see that 
\begin{align}
\frac{N_i(M)}{M} \tend{\rm a.s.}{M \to \infty } \frac{1}{\mu_{ii}} 
, \quad 
\frac{N(M)}{M} \tend{\rm a.s.}{M \to \infty } \frac{1}{\bar{m}} , 
\label{eq: NM/M}
\end{align}
and hence we obtain 
\begin{align}
A_{(i,j,k)}(M) \tend{\rm a.s.}{M \to \infty } a_{(i,j,k)} \frac{\mu_{ii}}{\sqrt{\bar{m}}} . 
\label{}
\end{align}
By Slutsky's theorem, the convergence \eqref{eq: vatildevq} 
for all constant $ \va = (a_{e})_{e \in \tilde{E}} \in \bR_0^{\tilde{E}} $ 
is reduced to the convergence 
\begin{align}
\vb \cdot \vY(M) 
\tend{\rm d}{M \to \infty } 
\cN(0,\vb \cdot V^{\vY} \vb) 
\label{eq: vbY}
\end{align}
for all constant $ \vb = (b_e)_{e \in \tilde{E}} \in \bR_0^{\tilde{E}} $, 
where, by \eqref{eq: barm}, 
\begin{align}
V^{\vY}_{ee'} 
:= \frac{\bar{m}}{\mu_{ii} \mu_{i'i'}} V^{\vq}_{ee'} 
= \nu_{i} 1_{\{ i=i' \}} q_{ij}(k) \rbra{ 1_{\{ j=j', \ k=k' \}} - q_{i'j'}(k') } . 
\label{}
\end{align}
for $ e = (i,j,k) $ and $ e' = (i',j',k') $. 
By the fact $ N(M)/M \tend{\rm a.s.}{M \to \infty } 1/\bar{m} > 0 $, 
the convergence \eqref{eq: vbY} is reduced to the convergence 
\begin{align}
Z(M) :=\frac{1}{\sqrt{M}} \sum_{n=1}^{M} \xi_n 
\tend{\rm d}{M \to \infty } & \cN(0,\vb \cdot V^{\vY} \vb) , 
\label{eq: vbZ}
\end{align}
where 
\begin{align}
\xi_n := \sum_{e = (i,j,k) \in \tilde{E}} b_e 
1_{\{J_{n-1}=i\}} \rbra{ 1_{\{ J_{n}=j, \ X_{n}=k \}} - q_{ij}(k) } . 
\label{}
\end{align}

Let us apply the martingale CLT (see, e.g., Theorem 8.2.8 of \cite{D}) to prove \eqref{eq: vbZ}. 
Set $ \cF_n = \sigma (J_l,X_l; l \le n ) $ for $ n \in \bN $. 
Then we see that 
$ \{ \xi_n \}_{n=1}^{\infty } $ 
is a martingale difference with respect to the filtration $ \{ \cF_n \}_{n=0}^{\infty } $, 
since 
\begin{align}
\bE_{i_0}[\xi_n \mid \cF_{n-1}] 
= \sum_{e = (i,j,k) \in \tilde{E}} b_e 1_{\{J_{n-1}=i\}} 
\bE_{i_0} \sbra{ \left. 1_{\{ J_{n}=j, \ X_{n}=k \}} - q_{ij}(k) \right| \cF_{n-1} } 
= 0 . 
\label{}
\end{align}
Since $ |\xi_n| \le \sum_e |b_e| < \infty $, 
it is obvious that 
\begin{align}
\frac{1}{M} \sum_{n=1}^M \bE_{i_0} \sbra{ \xi_n^2 \ ; \ |\xi_n| > \eps \sqrt{M} } 
\tend{}{M \to \infty } 0 
\quad \text{for all $ \eps>0 $}. 
\label{}
\end{align}
We have now only to prove that 
\begin{align}
F(M) := \frac{1}{M} \sum_{n=1}^M \bE_{i_0} \sbra{ \xi_n^2 \mid \cF_{n-1} } 
\tend{\rm a.s.}{M \to \infty } \vb \cdot V^{\vY} \vb . 
\label{}
\end{align}
Expanding the square, we obtain 
\begin{align}
F(M) 
= \frac{1}{M} \sum_{n=1}^M \sum_{\begin{subarray}{c} e = (i,j,k) \in \tilde{E} \\ e' = (i',j',k') \in \tilde{E} \end{subarray}} b_e b_{e'} 1_{\{J_{n-1}=i\}} 1_{\{ i=i' \}} E_{ee'}(n) 
\label{}
\end{align}
where, for $ i=i' $, 
\begin{align}
E_{ee'}(n) 
:=& \bE_{i_0} \sbra{ \rbra{ 1_{\{ J_{n}=j, \ X_{n}=k \}} - q_{ij}(k) } \rbra{ 1_{\{ J_{n}=j', \ X_{n}=k' \}} - q_{i'j'}(k') } \Bigm| \cF_{n-1} } 
\label{} \\
=& \bE_{i_0} \sbra{ 1_{\{ J_{n}=j, \ X_{n}=k \}} 1_{\{ J_{n}=j', \ X_{n}=k' \}} \Bigm| \cF_{n-1} } 
- q_{ij}(k) q_{i'j'}(k') 
\label{} \\
=& \bE_{i_0} \sbra{ 1_{\{ J_{n}=j, \ X_{n}=k \}} \Bigm| \cF_{n-1} } 1_{\{ j=j', \ k=k' \}} 
- q_{ij}(k) q_{i'j'}(k') 
\label{} \\
=& q_{ij}(k) \rbra{ 1_{\{ j=j', \ k=k' \}} - q_{i'j'}(k') } . 
\label{}
\end{align}
Using the fact from the ergodic theorem (see, e.g., Theorem 1.10.2 of \cite{N}) that 
\begin{align}
\frac{1}{M} \sum_{n=1}^M 1_{\{J_{n-1}=i \}} \tend{\rm a.s.}{M \to \infty } \nu_i , 
\label{}
\end{align}
we obtain 
\begin{align}
F(M) 
=& \sum_{\begin{subarray}{c} e = (i,j,k) \in \tilde{E} \\ e' = (i',j',k') \in \tilde{E} \end{subarray}} b_e b_{e'} \rbra{ \frac{1}{M} \sum_{n=1}^M 1_{\{J_{n-1}=i\}} } 1_{\{ i=i' \}} 
q_{ij}(k) \rbra{ 1_{\{ j=j', \ k=k' \}} - q_{i'j'}(k') } 
\label{} \\
\tend{\rm a.s.}{M \to \infty }& 
\sum_{\begin{subarray}{c} e = (i,j,k) \in \tilde{E} \\ e' = (i',j',k') \in \tilde{E} \end{subarray}} b_e b_{e'} \nu_i q_{ij}(k) \rbra{ 1_{\{ j=j', \ k=k' \}} - q_{i'j'}(k') } 
= \vb \cdot V^{\vY} \vb . 
\label{}
\end{align}
The proof is now complete. 
}

\section{Estimation of the matrix convolution inverse} \label{sec: MCI}

For two matrix sequences 
$ \vA = (\vA(k))_{k \in \bN} = (A_{ij}(k))_{(i,j,k) \in \tilde{E}} $ 
and 
$ \vB = (\vB(k))_{k \in \bN} = (B_{ij}(k))_{(i,j,k) \in \tilde{E}} $, 
we define the \emph{matrix convolution} of $ \vA $ and $ \vB $ 
by $ \vA * \vB = ((\vA * \vB)(k))_{k \in \bN} 
= ((\vA * \vB)_{ij}(k))_{(i,j,k) \in \tilde{E}} $ with 
\begin{align}
(\vA * \vB)(k) := \sum_{l=0}^k \vA(l) \vB(k-l) , 
\label{}
\end{align}
or equivalently, 
\begin{align}
(\vA * \vB)_{ij}(k) := \sum_{l=0}^k \sum_{u \in E} A_{iu}(l) B_{uj}(k-l) . 
\label{}
\end{align}
Define $ \vdI = (\vdI(k))_{k \in \bN} $ by 
\begin{align}
\vdI(k) = 
\begin{cases}
\vI & \text{for $ k = 0 $} \\
\vzero & \text{for $ k \ge 1 $} , 
\end{cases}
\label{}
\end{align}
where $ \vI $ and $ \vzero $ stand for the identity and zero matrices, respectively, 
so that $ \vdI $ is an identity of the matrix convolution: 
\begin{align}
\vdI * \vA = \vA * \vdI = \vA . 
\label{}
\end{align}
We denote the $ n $-fold matrix convolution of $ \vA $ by $ \vA^{(n)} $, 
that is, 
$ \vA^{(0)} = \vdI $ and $ \vA^{(n)} = \vA * \vA^{(n-1)} $ for $ n \ge 1 $. 
For any matrix-valued sequence $ \vA $ such that $ \vA^{(n)}(k) = \vzero $ whenever $ n>k $, 
we may define 
\begin{align}
(\vdI-\vA)^{(-1)} := \sum_{n=0}^{\infty } \vA^{(n)} , 
\label{}
\end{align}
which is the inverse of $ \vdI-\vA $ with respect to the matrix convolution: 
\begin{align}
(\vdI-\vA)^{(-1)} * (\vdI-\vA) = (\vdI-\vA) * (\vdI-\vA)^{(-1)} = \vdI 
\label{}
\end{align}
(see Proposition 3.4 of \cite{BL}). 

Note that the $ n $-fold matrix convolution of the transition probabilities, 
which we denote by 
$ \vq^{(n)} = (\vq^{(n)}(k))_{k \in \bN} = (q^{(n)}_{ij}(k))_{(i,j,k) \in \tilde{E}} $, 
satisfies 
\begin{align}
\bP_i(J_n=j, \ S_n=k) = q^{(n)}_{ij}(k) . 
\label{}
\end{align}
By the assumption $ \vq(0) = \vzero $, 
it is obvious that $ \vq^{(n)}(k) = \vzero $ whenever $ n>k $, 
and hence the matrix convolution inverse 
\begin{align}
\vpsi := (\vdI - \vq)^{(-1)} = \sum_{n=0}^{\infty } \vq^{(n)} 
\label{}
\end{align}
is well-defined. Since $ (\vdI - \vq) * \vpsi = \vpsi * (\vdI - \vq) = \vdI $, 
we have the renewal equation 
\begin{align}
\vpsi = \vdI + \vq * \vpsi = \vdI + \vpsi * \vq . 
\label{eq: RE vpsi}
\end{align}

For the estimator $ \hat{\vq}(M) $ introduced in \eqref{eq: hatvq} 
and its $ n $-fold matrix convolution $ \hat{\vq}^{(n)}(M) $, 
the inverse of the matrix convolution 
\begin{align}
\hat{\vpsi}(M) := (\vdI - \hat{\vq}(M))^{(-1)} = \sum_{n=0}^{\infty } \hat{\vq}^{(n)}(M) 
\label{}
\end{align}
is well-defined and satisfies the renewal equation 
\begin{align}
\hat{\vpsi}(M) 
= \vdI + \hat{\vq}(M) * \hat{\vpsi}(M) 
= \vdI + \hat{\vpsi}(M) * \hat{\vq}(M) . 
\label{eq: RE hatvpsi}
\end{align}
The strong consistency of $ \hat{\vpsi}(M) $ is known as follows. 

\begin{Thm}[Barbu--Limnios {\cite[Theorem 4.4]{BL}}]
Suppose {\bf (A)} is satisfied. Then 
\begin{align}
\hat{\vpsi}(M) \tend{\rm a.s.}{M \to \infty } \vpsi , 
\label{eq: SC vpsi}
\end{align}
under $ \bP_{i_0} $ for all initial value $ i_0 \in E $. 
\end{Thm}

The asymptotic normality of $ \hat{\vpsi}(M) = (\psi_{ij}(k,M))_{(i,j,k) \in \tilde{E}} $ 
in one-dimension has been obtained as follows. 

\begin{Thm}[Barbu--Limnios {\cite[Theorem 4.5]{BL}}] \label{1D vpsi}
Suppose {\bf (A)} is satisfied. 
Then, for any fixed $ (i,j,k) \in \tilde{E} $, it holds that 
\begin{align}
\sqrt{M} \rbra{ \hat{\psi}_{ij}(k,M) -\psi_{ij}(k) \Big. } 
\tend{\rm d}{M \to \infty } 
\cN(0,v^{\vpsi}(i,j,k)) 
\label{}
\end{align}
under $ \bP_{i_0} $ for all initial value $ i_0 \in E $, 
where the asymptotic variance $ v^{\vpsi}(i,j,k) $ is given as 
\begin{align}
v^{\vpsi}(i,j,k) 
:= \sum_{\tilde{i} \in E} \mu_{\tilde{i}\tilde{i}} \cbra{ \sum_{\tilde{j} \in E} \sbra{ (\psi_{i \, \tilde{i}} * \psi_{\tilde{j}j})^2 * q_{\tilde{i}\tilde{j}} }(k) 
- \sbra{ \sum_{\tilde{j} \in E} (\psi_{i \, \tilde{i}} * q_{\tilde{i}\tilde{j}} * \psi_{\tilde{j}j})(k) }^2 } . 
\label{eq: vvpsi}
\end{align}
\end{Thm}

Let us generalize Theorem \ref{1D vpsi} in multi-dimension, 
as a corollary of our Theorem \ref{thm: q}. 
Define the matrix sequences $ \bar{\vq} = (\bar{q}_{ij}(k))_{(i,j,k) \in \tilde{E}} $ 
and $ \vq^{\tilde{i}} = (q^{\tilde{i}}_{ij}(k))_{(i,j,k) \in \tilde{E}} $ 
for $ \tilde{i} \in E $ as 
\begin{align}
\bar{q}_{ij}(k) = \mu_{ii} q_{ij}(k) 
, \quad 
q^{\tilde{i}}_{ij}(k) = 1_{\{ i=\tilde{i} \}} q_{ij}(k) , 
\label{}
\end{align}
and define the matrix sequence $ \vdelta^{\tilde{e}} = (\delta^{\tilde{e}}_e)_{e \in \tilde{E}} $ 
for $ \tilde{e} \in \tilde{E} $ as 
\begin{align}
\delta^{\tilde{e}}_e = 1_{\{ e=\tilde{e} \}} 
\label{}
\end{align}
for $ e = (i,j,k) \in \tilde{E} $. 
For two matrix sequences $ \vA $ and $ \vB $, 
we define their tensor product $ \vA \otimes \vB $ 
as the $ \tilde{E} \times \tilde{E} $-matrix $ (A_e B_{e'})_{e,e' \in \tilde{E}} $. 

\begin{Thm} \label{thm: vpsi}
Suppose {\bf (A)} is satisfied. 
Then it holds that 
\begin{align}
\sqrt{M} \rbra{ \hat{\vpsi}(M) - \vpsi \Big. } 
\tend{\rm d}{M \to \infty } 
\cN(\vzero,V^{\vpsi}) 
\label{}
\end{align}
under $ \bP_{i_0} $ for all initial value $ i_0 \in E $, 
where the asymptotic variance matrix $ V^{\vpsi} = (V^{\vpsi}_{e e'})_{e,e' \in \tilde{E}} $ 
is given as 
\begin{align}
V^{\vpsi} 
= \sum_{\tilde{e} \in \tilde{E}} \bar{q}_{\tilde{e}} 
\rbra{ \vpsi * \vdelta^{\tilde{e}} * \vpsi }^{\otimes 2} 
- \sum_{\tilde{i} \in E} \mu_{\tilde{i} \tilde{i}} 
\rbra{ \vpsi * \vq^{\tilde{i}} * \vpsi }^{\otimes 2} . 
\label{eq: Vvpsi}
\end{align}
\end{Thm}

\Proof{
In what follows, we sometimes omit writing $ M $ from estimators, 
say $ \hat{\vq} $ for $ \hat{\vq}(M) $. 

By the renewal equations \eqref{eq: RE vpsi} and \eqref{eq: RE hatvpsi}, we have 
\begin{align}
\vdI = \rbra{ \vdI - \vq \Big. } * \vpsi = \hat{\vpsi} * \rbra{ \vdI - \hat{\vq} \Big. } , 
\label{}
\end{align}
and hence we have 
\begin{align}
\hat{\vpsi} - \vpsi 
=& \hat{\vpsi} * \vdI - \vdI * \vpsi 
\label{} \\
=& \hat{\vpsi} * \rbra{ \vdI - \vq \Big. } * \vpsi 
- \hat{\vpsi} * \rbra{ \vdI - \hat{\vq} \Big. } * \vpsi 
\label{} \\
=& \hat{\vpsi} * \rbra{ \hat{\vq} - \vq \Big. } * \vpsi , 
\label{}
\end{align}
which leads to the identity 
\begin{align}
\sqrt{M} \rbra{ \hat{\vpsi} - \vpsi } 
= \vpsi * \cbra{ \sqrt{M} \rbra{ \hat{\vq} - \vq \Big. } } * \vpsi 
+ \rbra{ \hat{\vpsi} - \vpsi \Big. } * \cbra{ \sqrt{M} \rbra{ \hat{\vq} - \vq \Big. } } * \vpsi . 
\label{eq: hatvpsi id}
\end{align}
By the strong consistency \eqref{eq: SC vpsi}, by the asymptotic normality \eqref{eq: q} 
and by Slutsky's theorem, we have 
\begin{align}
\rbra{ \hat{\vpsi} - \vpsi \Big. } 
* \cbra{ \sqrt{M} \rbra{ \hat{\vq} - \vq \Big. } } * \vpsi 
\tend{\rm d}{M \to \infty } \vzero . 
\label{eq: hatvpsi vzero}
\end{align}
Hence, by the asymptotic normality \eqref{eq: q} and by Slutsky's theorem, we obtain 
\begin{align}
\sqrt{M} \rbra{ \hat{\vpsi} - \vpsi \Big. } 
\tend{\rm d}{M \to \infty } 
\vpsi * \vZ^{\vq} * \vpsi , 
\label{}
\end{align}
where $ \vZ^{\vq} \dist \cN(\vzero,V^{\vq}) $. 
Since $ \vZ^{\vq} $ is a centered Gaussian system, 
we see that so is $ \vZ^{\vpsi} := \vpsi * \vZ^{\vq} * \vpsi $. 
The asymptotic variance matrix $ V^{\vpsi} $ of $ \vZ^{\vpsi} $ is given as 
\begin{align}
V^{\vpsi} 
=& \bE 
\sbra{ \rbra{ \vpsi * \vZ^{\vq} * \vpsi } \otimes \rbra{ \vpsi * \vZ^{\vq} * \vpsi } \Big. } 
\label{} \\
=& \sum_{\tilde{e},\tilde{e}' \in \tilde{E}} \bE \sbra{ Z^{\vq}_{\tilde{e}} Z^{\vq}_{\tilde{e}'} } 
\rbra{ \vpsi * \vdelta^{\tilde{e}} * \vpsi } \otimes \rbra{ \vpsi * \vdelta^{\tilde{e}'} * \vpsi } 
\label{} \\
=& \sum_{\tilde{e},\tilde{e}' \in \tilde{E}} V^{\vq}_{\tilde{e} \tilde{e}'} 
\rbra{ \vpsi * \vdelta^{\tilde{e}} * \vpsi } \otimes \rbra{ \vpsi * \vdelta^{\tilde{e}'} * \vpsi } 
\label{} \\
=& \sum_{\tilde{e} \in \tilde{E}} \bar{q}_{\tilde{e}} 
\rbra{ \vpsi * \vdelta^{\tilde{e}} * \vpsi }^{\otimes 2} 
- \sum_{\tilde{i} \in E} \mu_{\tilde{i} \tilde{i}} 
\rbra{ \vpsi * \vq^{\tilde{i}} * \vpsi }^{\otimes 2} . 
\label{}
\end{align}
The proof is complete. 
}

\section{Estimation of the distribution matrix sequence} \label{sec: DMS}

Let $ \cM $ denote the set of all random matrix sequences 
$ \vA = (A_{ij}(k))_{(i,j,k) \in \tilde{E}} $ and set 
\begin{align}
\cD = \cbra{ \vA \in \cM : 
\lim_{K \to \infty } \sum_{j \in E} \sum_{k \le K} A_{ij}(k) 
\ \text{exists in $ L^2 $ for all $ i \in E $} } . 
\label{}
\end{align}
We define the linear operator $ \cS : \cD \to \cM $ by 
\begin{align}
(\cS \vA)_{ij}(k) 
:= \lim_{K \to \infty } 1_{\{ i=j \}} \sum_{j' \in E} \sum_{k < k' \le K} A_{ij'}(k') 
\quad \text{in $ L^2 $} 
\label{}
\end{align}
with $ \cS \vA = ((\cS \vA)_{ij}(k))_{(i,j,k) \in \tilde{E}} $. 

For the semi-Markov chain $ Z_k = J_{N(k)} $ 
with $ N(k) = \max \{ n \in \bN : S_n \le k \} $, its distribution matrix sequence 
$ \vP = (P_{ij}(k))_{(i,j,k) \in \tilde{E}} $ has been given as 
\begin{align}
P_{ij}(k) := \bP_i(Z_k=j) 
= \bP_i(Z_k=j, \ X_1>k) + \bP_i(Z_k=j, \ X_1 \le k) . 
\label{}
\end{align}
Note that 
\begin{align}
\bP_i(Z_k=j, \ X_1>k) 
= 1_{\{ i=j \}} \bP_i(X_1>k) 
= 1_{\{ i=j \}} \sum_{j' \in E} \sum_{k'>k} q_{ij'}(k') 
= (\cS \vq)_{ij}(k) 
\label{}
\end{align}
and 
\begin{align}
\bP_i(Z_k=j, \ X_1 \le k) 
=& \sum_{j' \in E} \sum_{k' \le k} \bP_i(Z_k=j, \ J_1 = j' , \ X_1 = k') 
\label{} \\
=& \sum_{j' \in E} \sum_{k' \le k} q_{ij'}(k') \bP_{j'}(Z_{k-k'}=j) 
\label{} \\
=& \sum_{j' \in E} \sum_{k' \le k} q_{ij'}(k') P_{j'j}(k-k') 
= (\vq * \vP)_{ij}(k) . 
\label{}
\end{align}
We thus obtain the renewal equation 
\begin{align}
\vP = \cS \vq + \vq * \vP , 
\label{}
\end{align}
which can be solved as 
\begin{align}
\vP = (\vdI - \vq)^{(-1)} * \cS \vq = \vpsi * \cS \vq . 
\label{}
\end{align}

Let us introduce the estimator of the distribution matrix sequence by 
\begin{align}
\hat{\vP}(M) := \hat{\vpsi}(M) * \cS \hat{\vq}(M) . 
\label{}
\end{align}
The strong consistency of $ \hat{\vP}(M) $ is known as follows. 

\begin{Thm}[Barbu--Limnios {\cite[Theorem 4.6]{BL}}]
Suppose {\bf (A)} is satisfied. Then 
\begin{align}
\hat{\vP}(M) \tend{\rm a.s.}{M \to \infty } \vP , 
\label{eq: SC vP}
\end{align}
under $ \bP_{i_0} $ for all initial value $ i_0 \in E $. 
\end{Thm}

The asymptotic normality of $ \hat{\vP}(M) = (P_{ij}(k,M))_{(i,j,k) \in \tilde{E}} $ 
in one-dimension has been obtained as follows. 

\begin{Thm}[Barbu--Limnios {\cite[Theorem 4.7]{BL}}] \label{1D vP}
Suppose {\bf (A)} is satisfied. 
Then, for any fixed $ (i,j,k) \in \tilde{E} $, it holds that 
\begin{align}
\sqrt{M} \rbra{ \hat{P}_{ij}(k,M) -P_{ij}(k) \Big. } 
\tend{\rm d}{M \to \infty } 
\cN(0,v^{\vP}(i,j,k)) 
\label{}
\end{align}
under $ \bP_{i_0} $ for all initial value $ i_0 \in E $, 
where the asymptotic variance $ v^{\vP}(i,j,k) $ is given as 
\begin{align}
\begin{split}
v^{\vP}(i,j,k) 
:= \sum_{\tilde{i} \in E} \mu_{\tilde{i} \tilde{i}} \, \Biggl\{ & 
\sum_{\tilde{j} \in E} 
\sbra{ \rbra{ C^{ij}_{\tilde{i}\tilde{j}} - 1_{\{ \tilde{i}=j \}} \Psi_{ij} }^2 
* q_{\tilde{i} \tilde{j}} } (k) 
\\
&- \biggl[ 
\sum_{\tilde{j} \in E} C^{ij}_{\tilde{i}\tilde{j}} * q_{\tilde{i} \tilde{j} } 
- 1_{\{ \tilde{i}=j \}} \psi_{ij} * H_j 
\biggr]^2 (k) \Biggr\} , 
\end{split}
\label{eq: vvP}
\end{align}
where 
\begin{align}
\Psi_{ij}(k) = \sum_{k' \le k} \psi_{ij}(k') 
, \quad 
H_j(k) = \sum_{j' \in E} \sum_{k' \le k} q_{jj'}(k') 
\label{}
\end{align}
and 
\begin{align}
C^{ij}_{\tilde{i}\tilde{j}} = \psi_{i \, \tilde{i}} * \psi_{\tilde{j} j} * (1-H_j) . 
\label{}
\end{align}
\end{Thm}

Let us generalize Theorem \ref{1D vP} in multi-dimension, 
as a corollary of our Theorem \ref{thm: q}. 

\begin{Thm} \label{thm: vP}
Suppose {\bf (A)} is satisfied. 
Then it holds that 
\begin{align}
\sqrt{M} \rbra{ \hat{\vP}(M) - \vP \Big. } 
\tend{\rm d}{M \to \infty } 
\cN(\vzero,V^{\vP}) 
\label{}
\end{align}
under $ \bP_{i_0} $ for all initial value $ i_0 \in E $, 
where the asymptotic variance matrix $ V^{\vP} = (V^{\vP}_{e e'})_{e,e' \in \tilde{E}} $ 
is given as 
\begin{align}
\begin{split}
V^{\vP} 
=& \sum_{\tilde{e} \in \tilde{E}} \bar{q}_{\tilde{e}} 
\rbra{ \vpsi * \vdelta^{\tilde{e}} * \vpsi * \cS \vq 
+ \vpsi * \cS \vdelta^{\tilde{e}} }^{\otimes 2} 
\\
&- \sum_{\tilde{i} \in E} \mu_{\tilde{i} \tilde{i}} 
\rbra{ \vpsi * \vq^{\tilde{i}} * \vpsi * \cS \vq 
+ \vpsi * \cS \vq^{\tilde{i}} }^{\otimes 2} . 
\end{split}
\label{eq: VvP}
\end{align}
\end{Thm}

We need the following lemma. 

\begin{Lem} \label{lem: Z}
Let $ \vZ^{\vq} = (Z^{\vq}_e)_{e \in \tilde{E}} = (Z^{\vq}_{ij}(k))_{(i,j,k) \in \tilde{E}} 
\dist \cN(\vzero,V^{\vq}) $. Then 
\begin{align}
\lim_{K \to \infty } \sum_{j \in E} \sum_{k \le K} Z^{\vq}_{ij}(k) 
= 0 
\quad \text{in $ L^2 $} 
\label{eq: lemZ1}
\end{align}
for all $ i \in E $. Consequently, it holds that $ \vZ^{\vq} \in \cD $ and 
\begin{align}
(\cS \vZ^{\vq})_{ij}(k) = - 1_{\{ i=j \}} \sum_{j' \in E} \sum_{k' \le k} Z^{\vq}_{ij'}(k') 
\label{eq: lemZ2}
\end{align}
for all $ (i,j,k) \in \tilde{E} $. 
\end{Lem}

\Proof{
Expanding the square, we obtain 
\begin{align}
\bE \rbra{ \sum_{j \in E} \sum_{k \le K} Z^{\vq}_{ij}(k) }^2 
=& \sum_{j \in E} \sum_{k \le K} \sum_{j' \in E} \sum_{k' \le K} V_{(i,j,k)(i,j',k')} 
\label{} \\
=& \mu_{ii} \sum_{j \in E} \sum_{k \le K} q_{ij}(k) \sum_{j' \in E} \sum_{k' \le K} 
\rbra{ 1_{\{ j=j' , \ k=k' \}} - q_{ij'}(k') } 
\label{} \\
=& \mu_{ii} \sum_{j \in E} \sum_{k \le K} q_{ij}(k) 
\rbra{ 1 - \sum_{j' \in E} \sum_{k' \le K} q_{ij'}(k') } 
\label{} \\
\le& \mu_{ii} \sum_{j \in E} \sum_{k \in \bN} q_{ij}(k) \sum_{j' \in E} \sum_{k' > K} q_{ij'}(k') 
\tend{}{K \to \infty } 0 
\label{}
\end{align}
by the dominated convergence theorem. 
}

We also need the following lemma. 

\begin{Lem}
Suppose {\bf (A)} is satisfied. 
Then it holds that 
\begin{align}
\rbra{ \sqrt{M}(\hat{\vq}-\vq) , \ \cS \sqrt{M}(\hat{\vq}-\vq) \Big. } 
\tend{\rm d}{M \to \infty } 
\rbra{ \vZ^{\vq} , \ \cS \vZ^{\vq} \Big. } 
\label{}
\end{align}
in the product topology. 
\end{Lem}

\Proof{
Set $ \tilde{\vq} := \sqrt{M} (\hat{\vq} - \vq) $, 
so that $ \tilde{\vq} \tend{\rm d}{M \to \infty } \vZ^{\vq} $ 
by the asymptotic normality \eqref{eq: q}. 
It suffices to show that 
\begin{align}
\va \cdot \tilde{\vq} + \vb \cdot \cS \tilde{\vq} 
\tend{\rm d}{M \to \infty } 
\va \cdot \vZ^{\vq} + \vb \cdot \cS \vZ^{\vq} 
\label{}
\end{align}
for all constants $ \va,\vb \in \bR_0^{\tilde{E}} $. 
Since 
\begin{align}
\sum_{j' \in E} \sum_{k' \in \bN} \tilde{q}_{ij'}(k') 
=& \sqrt{M} \sum_{j' \in E} \sum_{k' \in \bN} \hat{q}_{ij'}(k') 
- \sqrt{M} \sum_{j' \in E} \sum_{k' \in \bN} q_{ij'}(k') 
\label{} \\
=& \sqrt{M} \cdot 1 - \sqrt{M} \cdot 1 =  0 , 
\label{}
\end{align}
we obtain 
\begin{align}
\va \cdot \tilde{\vq} + \vb \cdot \cS \tilde{\vq} 
=& \sum_{e \in \tilde{E}} a_e \tilde{q}_e 
+ \sum_{e=(i,j,k) \in \tilde{E}} b_e 1_{\{ i=j \}} 
\sum_{j' \in E} \sum_{k' > k} \tilde{q}_{ij'}(k') 
\label{} \\
=& \sum_{e \in \tilde{E}} a_e \tilde{q}_e 
- \sum_{e=(i,j,k) \in \tilde{E}} b_e 1_{\{ i=j \}} 
\sum_{j' \in E} \sum_{k' \le k} \tilde{q}_{ij'}(k') 
\label{} \\
\tend{\rm d}{M \to \infty } & 
\sum_{e \in \tilde{E}} a_e Z^{\vq}_e 
- \sum_{e=(i,j,k) \in \tilde{E}} b_e 1_{\{ i=j \}} 
\sum_{j' \in E} \sum_{k' \le k} Z^{\vq}_{ij'}(k') 
\label{eq: limit lemma} \\
=& \va \cdot \vZ^{\vq} + \vb \cdot \cS \vZ^{\vq} , 
\label{eq: limit lemma2}
\end{align}
where \eqref{eq: limit lemma} is obtained by applying \eqref{eq: lemZ1} to the finite sum, 
and where \eqref{eq: limit lemma2} is obtained by \eqref{eq: lemZ2}. 
The proof is complete. 
}

Let us proceed to the proof of Theorem \ref{thm: vP}. 

\Proof[Proof of Theorem \ref{thm: vP}]{
By \eqref{eq: hatvpsi id} and \eqref{eq: hatvpsi vzero}, we have 
\begin{align}
\sqrt{M} (\hat{\vpsi} - \vpsi) 
=& \vpsi * \sqrt{M} (\hat{\vq} - \vq) * \vpsi + o(1) , 
\label{}
\end{align}
where and in the sequel 
$ o(1) $ stands for an error term which vanishes in distribution as $ M \to \infty $. 
Hence we have 
\begin{align}
\sqrt{M} (\hat{\vP} - \vP) 
=& \sqrt{M} ( \hat{\vpsi} * \cS \hat{\vq} - \vpsi * \cS \vq ) 
\label{} \\
=& \sqrt{M} ( \hat{\vpsi} - \vpsi) * \cS \hat{\vq} + \sqrt{M} \vpsi * \cS ( \hat{\vq} - \vq ) 
\label{} \\
=& \sqrt{M} ( \hat{\vpsi} - \vpsi) * \cS \vq + o(1) 
+ \vpsi * \cS \sqrt{M}( \hat{\vq} - \vq ) 
\label{} \\
=& \vpsi * \sqrt{M} ( \hat{\vq} - \vq) * \vpsi * \cS \vq 
+ \vpsi * \cS \sqrt{M} ( \hat{\vq} - \vq ) + o(1) . 
\label{}
\end{align}
By the asymptotic normality \eqref{eq: q} and by Slutsky's theorem, we obtain 
\begin{align}
\sqrt{M} (\hat{\vP} - \vP) 
\tend{\rm d}{M \to \infty } 
\vpsi * \vZ^{\vq} * \vpsi * \cS \vq 
+ \vpsi * \cS \vZ^{\vq} . 
\label{}
\end{align}
The remainder of the proof is the same as that of Theorem \ref{thm: vpsi}. 
}

\section{Estimation of the reliability vector sequence} \label{sec: RVS}

Recall that the partition $ E = U \cup D $ 
and the reliability vector sequence $ \vR $ has been introduced 
in the end of Introduction. 

Let $ \tilde{U} := U \times U \times \bN $ 
and let $ \vq_{UU} = ((q_{UU})_{ij}(k))_{(i,j,k) \in \tilde{U}} $ 
denote the restriction of $ \vq $ on $ \tilde{U} $, i.e., 
\begin{align}
(q_{UU})_{ij}(k) := q_{ij}(k) 
\quad \text{for $ (i,j,k) \in \tilde{U} $.} 
\label{}
\end{align}
Let $ \vq^{(n)}_{UU} = ((q^{(n)}_{UU})_{ij}(k))_{(i,j,k) \in \tilde{U}} $ 
denote the $ n $-fold matrix convolution of $ \vq_{UU} $, which satisfies 
\begin{align}
\bP_i(J_n=j, \ S_n=k , \ T_D > k) = (q^{(n)}_{UU})_{ij}(k) . 
\label{}
\end{align}
We denote the matrix convolution inverse of $ \vdI - \vq_{UU} $ by 
\begin{align}
\vpsi_{UU} := (\vdI - \vq_{UU})^{(-1)} = \sum_{n=0}^{\infty } \vq_{UU}^{(n)} . 
\label{}
\end{align}
Note that $ \vq^{(n)}_{UU} $ for $ n \ge 2 $ and $ \vpsi_{UU} $ 
\emph{differ} from the restriction of $ \vq^{(n)} $ and $ \vpsi $ on $ \tilde{U} $, respectively. 

The reliability vector sequence $ \vR = (R_i(k))_{(i,k) \in U \times \bN} $ 
has been given as 
\begin{align}
R_i(k) := \bP_i(T_D>k) 
= \bP_i(T_D>k , \ X_1>k) + \bP_i(T_D>k , \ X_1 \le k) . 
\label{}
\end{align}
Note that 
\begin{align}
\bP_i(T_D>k , \ X_1>k) = \bP_i(X_1>k) = (\cS \vq)_{ii}(k) 
= ((\cS \vq)_U)_i(k) , 
\label{}
\end{align}
where we write $ \vA_U = ((\vA_U)_i(k))_{(i,k) \in U \times \bN} 
:= (A_{ii}(k))_{(i,k) \in U \times \bN} $ 
for a matrix sequence $ \vA = (A_{ij}(k))_{(i,j,k) \in \tilde{E}} $, 
and 
\begin{align}
\bP_i(T_D>k , \ X_1 \le k) 
=& \sum_{j' \in U} \sum_{k' \le k} \bP_i(T_D>k , \ J_1 = j' , \ X_1 = k') 
\label{} \\
=& \sum_{j' \in U} \sum_{k' \le k} q_{ij'}(k') \bP_{j'}(T_D>k-k') 
\label{} \\
=& \sum_{j' \in U} \sum_{k' \le k} q_{ij'}(k') R_{j'}(k-k') 
= (\vq_{UU} * \vR)_i(k) . 
\label{}
\end{align}
We thus obtain the renewal equation 
\begin{align}
\vR = (\cS \vq)_U + \vq_{UU} * \vR , 
\label{}
\end{align}
which can be solved as 
\begin{align}
\vR = (\vdI - \vq_{UU})^{(-1)} * (\cS \vq)_U = \vpsi_{UU} * (\cS \vq)_U . 
\label{}
\end{align}

Define 
$ \hat{\vq}_{UU} = ((\hat{q}_{UU})_{ij}(k))_{(i,j,k) \in \tilde{U}} $ 
by $ (\hat{q}_{UU})_{ij}(k) := \hat{q}_{ij}(k) $ for $ (i,j,k) \in \tilde{U} $ and 
\begin{align}
\hat{\vpsi}_{UU} := (\vdI - \hat{\vq}_{UU})^{(-1)} 
= \sum_{n=0}^{\infty } \hat{\vq}_{UU}^{(n)} , 
\label{}
\end{align}
where $ \hat{\vq}^{(n)}_{UU} = ((\hat{q}^{(n)}_{UU})_{ij}(k))_{(i,j,k) \in \tilde{U}} $ 
denote the $ n $-fold matrix convolution of $ \hat{\vq}_{UU} $. 
Let us introduce the estimator of the reliability vector sequence by 
\begin{align}
\hat{\vR} := \hat{\vpsi}_{UU} * (\cS \hat{\vq})_U . 
\label{}
\end{align}
The strong consistency of $ \hat{\vR} $ is known as follows. 

\begin{Thm}[Barbu--Limnios {\cite[Theorem 5.1]{BL}}]
Suppose {\bf (A)} is satisfied. Then 
\begin{align}
\hat{\vR} \tend{\rm a.s.}{M \to \infty } \vR , 
\label{eq: SC vR}
\end{align}
under $ \bP_{i_0} $ for all initial value $ i_0 \in U $. 
\end{Thm}

The asymptotic normality of $ \hat{\vR} = (\hat{R}_i(k,M))_{(i,k) \in U \times \bN} $ 
in one-dimension is as follows, 
which corrects the wrong statement of Barbu--Limnios \cite[Theorem 5.1]{BL}. 

\begin{Thm} \label{1D vR}
Suppose {\bf (A)} is satisfied. 
Then, for any fixed $ (i,k) \in U \times \bN $, it holds that 
\begin{align}
\sqrt{M} \rbra{ \hat{R}_i(k,M) -R_i(k) \Big. } 
\tend{\rm d}{M \to \infty } 
\cN(0,v^{\vR}(i,k)) 
\label{}
\end{align}
under $ \bP_{i_0} $ for all initial value $ i_0 \in U $, 
where the asymptotic variance $ v^{\vR}(i,k) $ is given as 
\begin{align}
\begin{split}
v^{\vR}(i,k) :=& \sum_{\tilde{i} \in U} \mu_{\tilde{i}\tilde{i}} 
\Biggl\{ \sum_{\tilde{j} \in U} \sbra{ \rbra{ D^i_{\tilde{i}\tilde{j}} - (\Psi_{UU})_{i \, \tilde{i}} }^2 * q_{\tilde{i}\tilde{j}} } (k) 
\\
&- \sum_{\tilde{j} \in U} \sbra{ D^i_{\tilde{i}\tilde{j}} * q_{\tilde{i}\tilde{j}} 
- (\psi_{UU})_{i \, \tilde{i}} * H_{\tilde{i}} }^2 (k) \Biggr\} 
\\
&+ \sum_{\tilde{i} \in U} \mu_{\tilde{i}\tilde{i}} 
\sum_{\tilde{j} \in D} \sbra{ \rbra{ (\Psi_{UU})_{i \, \tilde{i}}(k) - (\Psi_{UU})_{i \, \tilde{i}} }^2 * q_{\tilde{i} \tilde{j}} }(k) , 
\end{split}
\label{eq: vvR}
\end{align}
where 
\begin{align}
(\Psi_{UU})_{ij}(k) = \sum_{k' \le k} (\psi_{UU})_{ij}(k') 
\label{}
\end{align}
and 
\begin{align}
D^i_{\tilde{i}\tilde{j}} = \sum_{j \in U} (\psi_{UU})_{i \, \tilde{i}} * (\psi_{UU})_{\tilde{j} j} * (1-H_j) . 
\label{}
\end{align}
\end{Thm}

Let us generalize Theorem \ref{1D vR} in multi-dimension, 
as a corollary of our Theorem \ref{thm: q}. 

\begin{Thm} \label{thm: vR}
Suppose {\bf (A)} is satisfied. 
Then it holds that 
\begin{align}
\sqrt{M} \rbra{ \hat{\vR} - \vR \Big. } 
\tend{\rm d}{M \to \infty } 
\cN(\vzero,V^{\vR}) 
\label{eq: cdist vR}
\end{align}
under $ \bP_{i_0} $ for all initial value $ i_0 \in U $, 
where the asymptotic variance matrix $ V^{\vR} = (V^{\vR}_{e e'})_{e,e' \in U \times \bN} $ 
is given as 
\begin{align}
\begin{split}
V^{\vR} 
=& \sum_{\tilde{e} \in \tilde{U}} \bar{q}_{\tilde{e}} 
\rbra{ \vpsi_{UU} * \vdelta^{\tilde{e}} * \vpsi_{UU} * (\cS \vq)_U 
+ \vpsi_{UU} * (\cS \vdelta^{\tilde{e}})_U }^{\otimes 2} 
\\
&- \sum_{\tilde{i} \in U} \mu_{\tilde{i} \tilde{i}} 
\rbra{ \vpsi_{UU} * \vq^{\tilde{i}}_{UU} * \vpsi_{UU} * (\cS \vq)_U 
+ \vpsi_{UU} * (\cS \vq^{\tilde{i}})_U }^{\otimes 2} 
\\
&+ \sum_{\tilde{e} \in U \times D \times \bN} \bar{q}_{\tilde{e}} 
\rbra{ \vpsi_{UU} * (\cS \vdelta^{\tilde{e}})_U }^{\otimes 2} . 
\end{split}
\label{eq: VvR}
\end{align}
\end{Thm}

\Proof{
By the same arguments as \eqref{eq: hatvpsi id} and \eqref{eq: hatvpsi vzero}, we have 
\begin{align}
\sqrt{M} (\hat{\vpsi}_{UU} - \vpsi_{UU}) 
=& \vpsi_{UU} * \sqrt{M} (\hat{\vq}_{UU} - \vq_{UU}) * \vpsi_{UU} + o(1) . 
\label{}
\end{align}
Hence we have 
\begin{align}
\sqrt{M} (\hat{\vR} - \vR) 
=& \sqrt{M} ( \hat{\vpsi}_{UU} * (\cS \hat{\vq})_U - \vpsi_{UU} * (\cS \vq)_U ) 
\label{} \\
=& \sqrt{M} ( \hat{\vpsi}_{UU} - \vpsi_{UU}) * (\cS \hat{\vq})_U 
+ \sqrt{M} \vpsi_U * (\cS ( \hat{\vq} - \vq ))_U 
\label{} \\
=& \sqrt{M} ( \hat{\vpsi}_{UU} - \vpsi_{UU}) * (\cS \vq)_U + o(1) 
+ \vpsi_{UU} * (\cS \sqrt{M}( \hat{\vq} - \vq ))_U 
\label{} \\
\begin{split}
=& \vpsi_{UU} * \sqrt{M} ( \hat{\vq}_{UU} - \vq_{UU}) * \vpsi_{UU} * (\cS \vq)_U 
\\
&+ \vpsi_{UU} * (\cS \sqrt{M} ( \hat{\vq} - \vq ))_U + o(1) . 
\end{split}
\label{}
\end{align}
By the asymptotic normality \eqref{eq: q} and by Slutsky's theorem, we obtain 
\begin{align}
\sqrt{M} (\hat{\vR} - \vR) 
\tend{\rm d}{M \to \infty } 
\vpsi_{UU} * \vZ^{\vq}_{UU} * \vpsi_{UU} * (\cS \vq)_U 
+ \vpsi_{UU} * (\cS \vZ^{\vq})_U . 
\label{}
\end{align}
By the same arguments as the proof of Theorem \ref{thm: vpsi}, 
we obtain \eqref{eq: cdist vR} with the asymptotic covariance matrix being given by 
\begin{align}
\begin{split}
V^{\vR} 
=& \sum_{\tilde{e} \in \tilde{E}} \bar{q}_{\tilde{e}} 
\rbra{ \vpsi_{UU} * \vdelta^{\tilde{e}}_{UU} * \vpsi_{UU} * (\cS \vq)_U 
+ \vpsi_{UU} * (\cS \vdelta^{\tilde{e}})_U }^{\otimes 2} 
\\
&- \sum_{\tilde{i} \in E} \mu_{\tilde{i} \tilde{i}} 
\rbra{ \vpsi_{UU} * \vq^{\tilde{i}}_{UU} * \vpsi_{UU} * (\cS \vq)_U 
+ \vpsi_{UU} * (\cS \vq^{\tilde{i}})_U }^{\otimes 2} , 
\end{split}
\label{}
\end{align}
where $ \vdelta^{\tilde{e}}_{UU} $ stands for the restriction 
of $ \vdelta^{\tilde{e}} $ on $ \tilde{U} $. 
Since we easily see that 
\begin{align}
\vpsi_{UU} * \vdelta^{\tilde{e}}_{UU} * \vpsi_{UU} = \vzero 
\quad \text{for $ \tilde{e} \notin \tilde{U} $} , 
\label{}
\end{align}
\begin{align}
\vpsi_{UU} * (\cS \vdelta^{\tilde{e}})_U = \vzero 
\quad \text{for $ \tilde{e} \in D \times E \times \bN $} 
\label{}
\end{align}
and 
\begin{align}
\vpsi_{UU} * \vq^{\tilde{i}}_{UU} * \vpsi_{UU} 
= \vpsi_{UU} * (\cS \vq^{\tilde{i}})_U 
= \vzero 
\quad \text{for $ \tilde{i} \in D $} , 
\label{}
\end{align}
we obtain \eqref{eq: VvR}. 
}

\stepcounter{section}
\section*{Appendix: Computation of the asymptotic variances}

\subsection{The asymptotic variance for the matrix convolution inverse}

Let us check that $ V^{\vpsi}_{ee'} $ in \eqref{eq: Vvpsi} for $ e=e'=(i,j,k) \in \tilde{E} $ 
coincides with $ v^{\vpsi}(i,j,k) $ in \eqref{eq: vvpsi}. 
Set 
\begin{align}
V^1 := \sum_{\tilde{e} \in \tilde{E}} \bar{q}_{\tilde{e}} 
\rbra{ \vpsi * \vdelta^{\tilde{e}} * \vpsi }^{\otimes 2} 
, \quad 
V^2 := \sum_{\tilde{i} \in E} \mu_{\tilde{i} \tilde{i}} 
\rbra{ \vpsi * \vq^{\tilde{i}} * \vpsi }^{\otimes 2} , 
\label{}
\end{align}
so that $ V^{\vpsi} = V^1 - V^2 $. 
For $ V^1_{ee} $, we have 
\begin{align}
V^1_{ee} 
=& \sum_{\tilde{e} \in \tilde{E}} \bar{q}_{\tilde{e}} 
\sbra{ \rbra{ \vpsi * \vdelta^{\tilde{e}} * \vpsi }_e }^2 
= \sum_{\tilde{e} = (\tilde{i},\tilde{j},\tilde{k}) \in \tilde{E}} {q}_{\tilde{e}} 
\sbra{ (\psi_{i \, \tilde{i}} * \psi_{\tilde{j} j})(k-\tilde{k}) }^2 
\label{} \\
=& \sum_{\tilde{i} \in E } \mu_{\tilde{i}\tilde{i}} \sum_{\tilde{j} \in E } 
\sum_{\tilde{k} \in \bN } 
q_{\tilde{i} \tilde{j}}(\tilde{k}) 
\sbra{ (\psi_{i \, \tilde{i}} * \psi_{\tilde{j} j} )(k-\tilde{k}) }^2 
\label{} \\
=& \sum_{\tilde{i} \in E } \mu_{\tilde{i}\tilde{i}} \sum_{\tilde{j} \in E } 
\sbra{ (\psi_{i \, \tilde{i}} * \psi_{\tilde{j} j} )^2 * q_{\tilde{i} \tilde{j}} }(k) . 
\label{}
\end{align}
For $ V^2_{ee} $, we have 
\begin{align}
V^2_{ee} 
=& \sum_{\tilde{i} \in E} \mu_{\tilde{i}\tilde{i}} \sbra{ \rbra{ \vpsi * \vq^{\tilde{i}} * \vpsi }_e }^2 
= \sum_{\tilde{i} \in E} \mu_{\tilde{i}\tilde{i}} \sbra{ \sum_{\tilde{j} \in E} \rbra{ \psi_{i \, \tilde{i}} * q_{\tilde{i} \tilde{j}} * \psi_{\tilde{j} j} } (k) }^2 . 
\label{}
\end{align}
Hence we obtain $ V^{\vpsi}_{ee} = v^{\vpsi}(i,j,k) $.

\subsection{The asymptotic variance for the distribution matrix sequence}

Let us check that $ V^{\vP}_{ee'} $ in \eqref{eq: VvP} for $ e=e'=(i,j,k) \in \tilde{E} $ 
coincides with $ v^{\vP}(i,j,k) $ in \eqref{eq: vvP}. 
Set 
\begin{align}
V^1 :=& \sum_{\tilde{e} \in \tilde{E}} \bar{q}_{\tilde{e}} 
\rbra{ \vpsi * \vdelta^{\tilde{e}} * \vpsi * \cS \vq 
+ \vpsi * \cS \vdelta^{\tilde{e}} }^{\otimes 2} , 
\label{} \\
V^2 :=& \sum_{\tilde{i} \in E} \mu_{\tilde{i} \tilde{i}} 
\rbra{ \vpsi * \vq^{\tilde{i}} * \vpsi * \cS \vq 
+ \vpsi * \cS \vq^{\tilde{i}} }^{\otimes 2} , 
\label{}
\end{align}
so that $ V^{\vP} = V^1 - V^2 $. 
Since 
\begin{align}
\rbra{ \vpsi * \vdelta^{\tilde{e}} * \vpsi * \cS \vq }_e 
=& \sum_{k'' \in \bN: \ \tilde{k} \le k'' \le k} 
\rbra{ \psi_{i \, \tilde{i}} * \psi_{\tilde{j}j} }(k''-\tilde{k}) 
\sum_{j' \in E} \sum_{k'>k-k''} q_{jj'}(k') 
\label{} \\
=& \sum_{k'' \in \bN: \ \tilde{k} \le k'' \le k} 
\rbra{ \psi_{i \, \tilde{i}} * \psi_{\tilde{j}j} }(k''-\tilde{k}) 
\cdot (1 - H_j(k-k'')) 
\label{} \\
=& \rbra{ \psi_{i \, \tilde{i}} * \psi_{\tilde{j}j} * (1 - H_j) }(k-\tilde{k}) 
= C^{ij}_{\tilde{i}\tilde{j}}(k-\tilde{k}) 
\label{}
\end{align}
and 
\begin{align}
\rbra{ \vpsi * \cS \vdelta^{\tilde{e}} }_e 
=& 1_{\{ \tilde{i}=j \}} \sum_{k' \le k} \psi_{i \, \tilde{i}}(k') 1_{\{ k' > k-\tilde{k} \ge 0 \}} 
\label{} \\
=& 1_{\{ \tilde{i}=j \}} \rbra{ \Psi_{ij}(k) - \Psi_{ij}(k-\tilde{k}) } 1_{\{ \tilde{k} \le k \}} , 
\label{}
\end{align}
we have 
\begin{align}
V^1_{ee} 
=& \sum_{\tilde{i} \in E } \mu_{\tilde{i}\tilde{i}} \sum_{\tilde{j} \in E } 
\sum_{\tilde{k} \in \bN } 
q_{\tilde{i} \tilde{j}}(\tilde{k}) 
\sbra{ \rbra{ \vpsi * \vdelta^{\tilde{e}} * \vpsi * \cS \vq }_e
+ \rbra{ \vpsi * \cS \vdelta^{\tilde{e}} }_e }^2 
\label{} \\
=& \sum_{\tilde{i} \in E } \mu_{\tilde{i}\tilde{i}} \sum_{\tilde{j} \in E } 
\sbra{ \rbra{ C^{ij}_{\tilde{i}\tilde{j}} 
+ 1_{\{ \tilde{i}=j \}} \rbra{ \Psi_{ij}(k) - \Psi_{ij} } }^2 * q_{\tilde{i} \tilde{j}} } (k) . 
\label{}
\end{align}
Since 
\begin{align}
\rbra{ \vpsi * \vq^{\tilde{i}} * \vpsi * \cS \vq }_e 
=& \sum_{\tilde{j} \in E} \sum_{k'' \in \bN} 
\rbra{ \psi_{i \, \tilde{i}} * q_{\tilde{i}\tilde{j}} * \psi_{\tilde{j}j} }(k'') 
\sum_{j' \in E} \sum_{k'>k-k''} q_{jj'}(k') 
\label{} \\
=& \sum_{\tilde{j} \in E} \sbra{ \psi_{i \, \tilde{i}} * q_{\tilde{i}\tilde{j}} * \psi_{\tilde{j}j} * (1-H_j) }(k)
\label{} \\
=& \sum_{\tilde{j} \in E} \sbra{ C^{ij}_{\tilde{i}\tilde{j}} * q_{\tilde{i}\tilde{j}} }(k) 
\label{}
\end{align}
and 
\begin{align}
\rbra{ \vpsi * \cS \vq^{\tilde{i}} }_e 
=& 1_{\{ \tilde{i}=j \}} \sum_{k'' \le k} \psi_{i \, \tilde{i}}(k'')  
\sum_{j' \in E} \sum_{k'>k-k''} q_{jj'}(k') 
\label{} \\
=& 1_{\{ \tilde{i}=j \}} \sbra{ \psi_{ij} * (1-H_j) }(k) , 
\label{}
\end{align}
we have 
\begin{align}
V^2_{ee} 
=& \sum_{\tilde{i} \in E } \mu_{\tilde{i}\tilde{i}} 
\sbra{ \rbra{ \vpsi * \vq^{\tilde{i}} * \vpsi * \cS \vq }_e 
+ \rbra{ \vpsi * \cS \vq^{\tilde{i}} }_e }^2 , 
\label{} \\
=& \sum_{\tilde{i} \in E } \mu_{\tilde{i}\tilde{i}} 
\sbra{ \rbra{ \sum_{\tilde{j} \in E} C^{ij}_{\tilde{i}\tilde{j}} * q_{\tilde{i}\tilde{j}} 
+ 1_{\{ \tilde{i}=j \}} \psi_{ij} * (1-H_j) }^2 } (k) . 
\label{}
\end{align}
Since $ \psi_{ij} * 1 = \Psi_{ij} $, $ \sum_{\tilde{j} \in E} q_{i\tilde{j}} * 1 = 1 $ 
and $ \sum_{\tilde{j} \in E} \Psi_{ij} * q_{\tilde{i}\tilde{j}} = \psi_{ij} * H_{\tilde{i}} $, 
we have 
\begin{align}
\begin{split}
& \sum_{\tilde{j} \in E } 
\sbra{ \rbra{ C^{ij}_{\tilde{i}\tilde{j}} 
+ 1_{\{ \tilde{i}=j \}} \rbra{ \Psi_{ij}(k) - \Psi_{ij} } }^2 * q_{\tilde{i} \tilde{j}} } (k) 
\\
&- 
\sbra{ \rbra{ \sum_{\tilde{j} \in E} C^{ij}_{\tilde{i}\tilde{j}} * q_{\tilde{i}\tilde{j}} 
+ 1_{\{ \tilde{i}=j \}} \psi_{ij} * (1-H_j) }^2 } (k) 
\end{split}
\label{} \\
\begin{split}
=& \sum_{\tilde{j} \in E } 
\sbra{ \rbra{ C^{ij}_{\tilde{i}\tilde{j}} 
- 1_{\{ \tilde{i}=j \}} \Psi_{ij} }^2 * q_{\tilde{i} \tilde{j}} } (k) 
\\
&- 
\sbra{ \rbra{ \sum_{\tilde{j} \in E} C^{ij}_{\tilde{i}\tilde{j}} * q_{\tilde{i}\tilde{j}} 
- 1_{\{ \tilde{i}=j \}} \psi_{ij} * H_j }^2 } (k) . 
\end{split}
\label{}
\end{align}
Hence we obtain $ V^{\vP}_{ee} = v^{\vP}(i,j,k) $.

\subsection{The asymptotic variance for the reliability vector sequence}

Let us check that $ V^{\vR}_{ee'} $ in \eqref{eq: VvR} for $ e=e'=(i,k) \in U \times \bN $ 
coincides with $ v^{\vR}(i,k) $ in \eqref{eq: vvR}. 
Set 
\begin{align}
V^1 :=& \sum_{\tilde{e} \in \tilde{U}} \bar{q}_{\tilde{e}} 
\rbra{ \vpsi_{UU} * \vdelta^{\tilde{e}} * \vpsi_{UU} * (\cS \vq)_U 
+ \vpsi_{UU} * (\cS \vdelta^{\tilde{e}})_U }^{\otimes 2} , 
\label{} \\
V^2 :=& \sum_{\tilde{i} \in U} \mu_{\tilde{i} \tilde{i}} 
\rbra{ \vpsi_{UU} * \vq^{\tilde{i}}_{UU} * \vpsi_{UU} * (\cS \vq)_U 
+ \vpsi_{UU} * (\cS \vq^{\tilde{i}})_U }^{\otimes 2} , 
\label{} \\
V^3 :=& \sum_{\tilde{e} \in U \times D \times \bN} \bar{q}_{\tilde{e}} 
\rbra{ \vpsi_{UU} * (\cS \vdelta^{\tilde{e}})_U }^{\otimes 2} , 
\label{}
\end{align}
so that $ V^{\vR} = V^1 - V^2 + V^3 $. Since 
\begin{align}
& \rbra{ \vpsi_{UU} * \vdelta^{\tilde{e}} * \vpsi_{UU} * (\cS \vq)_U }_e 
\label{} \\
=& \sum_{j \in U} \sum_{k'' \in \bN : \ \tilde{k} \le k'' \le k} \rbra{ (\vpsi_{UU})_{i \, \tilde{i}} * (\vpsi_{UU})_{\tilde{j}j}}(k''-\tilde{k}) 
\sum_{j' \in E} \sum_{k'>k-k''} q_{jj'}(k') 
\label{} \\
=& \sum_{j \in U} \rbra{ (\vpsi_{UU})_{i \, \tilde{i}} * (\vpsi_{UU})_{\tilde{j}j} * (1-H_j) } (k-\tilde{k}) = D^{i}_{\tilde{i}\tilde{j}}(k-\tilde{k}) 
\label{}
\end{align}
and 
\begin{align}
\rbra{ \vpsi_{UU} * (\cS \vdelta^{\tilde{e}})_U }_e 
=& \sum_{k' \le k} (\psi_{UU})_{i \, \tilde{i}}(k') 1_{\{ k' > k-\tilde{k} \ge 0 \}} 
\label{} \\
=& \rbra{ (\Psi_{UU})_{i \, \tilde{i}}(k) - (\Psi_{UU})_{i \, \tilde{i}}(k-\tilde{k}) } 
1_{\{ \tilde{k} \le k \}} , 
\label{}
\end{align}
we have 
\begin{align}
V^1_{ee} 
=& \sum_{\tilde{i} \in U} \mu_{\tilde{i}\tilde{i}} 
\sum_{\tilde{j} \in U} \sum_{\tilde{k} \in \bN} 
q_{\tilde{i} \tilde{j}}(\tilde{k}) 
\sbra{ \rbra{ \vpsi_{UU} * \vdelta^{\tilde{e}} * \vpsi_{UU} * (\cS \vq)_U }_e 
+ \rbra{ \vpsi_{UU} * (\cS \vdelta^{\tilde{e}})_U }_e }^2 
\label{} \\
=& \sum_{\tilde{i} \in U} \mu_{\tilde{i}\tilde{i}} 
\sum_{\tilde{j} \in U} \sbra{ \rbra{ D^{i}_{\tilde{i}\tilde{j}} + (\Psi_{UU})_{i \, \tilde{i}}(k) - (\Psi_{UU})_{i \, \tilde{i}} }^2 * q_{\tilde{i} \tilde{j}} }(k) 
\label{}
\end{align}
and 
\begin{align}
V^3_{ee} 
=& \sum_{\tilde{i} \in U} \mu_{\tilde{i}\tilde{i}} 
\sum_{\tilde{j} \in D} \sum_{\tilde{k} \in \bN} 
q_{\tilde{i} \tilde{j}}(\tilde{k}) 
\sbra{ \rbra{ \vpsi_{UU} * (\cS \vdelta^{\tilde{e}})_U }_e }^2 , 
\label{} \\
=& \sum_{\tilde{i} \in U} \mu_{\tilde{i}\tilde{i}} 
\sum_{\tilde{j} \in D} \sbra{ \rbra{ (\Psi_{UU})_{i \, \tilde{i}}(k) - (\Psi_{UU})_{i \, \tilde{i}} }^2 * q_{\tilde{i} \tilde{j}} }(k) . 
\label{}
\end{align}
Since 
\begin{align}
& \rbra{ \vpsi_{UU} * \vq^{\tilde{i}}_{UU} * \vpsi_{UU} * (\cS \vq)_U }_e 
\label{} \\
=& \sum_{\tilde{j} \in U} \sum_{k'' \in \bN} \sum_{j \in U} \rbra{ (\psi_{UU})_{i \, \tilde{i}} * q_{\tilde{i} \tilde{j}} * (\psi_{UU})_{\tilde{j} j} }(k'') \sum_{j' \in E} \sum_{k'>k-k''} q_{jj'}(k') 
\label{} \\
=& \sum_{\tilde{j} \in U} \sbra{ (\psi_{UU})_{i \, \tilde{i}} * q_{\tilde{i} \tilde{j}} * (\psi_{UU})_{\tilde{j} j} * (1-H_j) }(k) 
\label{} \\
=& \sum_{\tilde{j} \in U} \sbra{ D^{i}_{\tilde{i} \tilde{j}} * q_{\tilde{i} \tilde{j}} }(k) 
\label{}
\end{align}
and 
\begin{align}
\rbra{ \vpsi_{UU} * (\cS \vq^{\tilde{i}})_U }_e 
=& \sum_{k'' \le k} (\psi_{UU})_{i \, \tilde{i}}(k'') \sum_{j' \in E} \sum_{k' > k-k''} q_{\tilde{i} j'}(k') 
\label{} \\
=& \sbra{ (\psi_{UU})_{i \, \tilde{i}} * (1-H_{\tilde{i}}) }(k) , 
\label{}
\end{align}
we have 
\begin{align}
V^2_{ee} 
=& \sum_{\tilde{i} \in U} \mu_{\tilde{i}\tilde{i}} 
\sbra{ \rbra{ \vpsi_{UU} * \vq^{\tilde{i}}_{UU} * \vpsi_{UU} * (\cS \vq)_U }_e + \rbra{ \vpsi_{UU} * (\cS \vq^{\tilde{i}})_U }_e }^2 
\label{} \\
=& \sum_{\tilde{i} \in U} \mu_{\tilde{i}\tilde{i}} 
\sbra{ \rbra{ \sum_{\tilde{j} \in U} D^{i}_{\tilde{i} \tilde{j}} * q_{\tilde{i} \tilde{j}} + (\psi_{UU})_{i \, \tilde{i}} * (1-H_{\tilde{i}}) }^2 }(k) . 
\label{}
\end{align}
Since $ (\psi_{UU})_{ij} * 1 = (\Psi_{UU})_{ij} $, 
$ \sum_{\tilde{j} \in E} q_{i\tilde{j}} * 1 = 1 $ 
and $ \sum_{\tilde{j} \in E} (\Psi_{UU})_{ij} * q_{\tilde{i}\tilde{j}} 
= (\psi_{UU})_{ij} * H_{\tilde{i}} $, 
we have 
\begin{align}
\begin{split}
& \sum_{\tilde{j} \in U } 
\sbra{ \rbra{ D^{i}_{\tilde{i}\tilde{j}} + (\Psi_{UU})_{i \, \tilde{i}}(k) - (\Psi_{UU})_{i \, \tilde{i}} }^2 * q_{\tilde{i} \tilde{j}} } (k) 
\\
&- 
\sbra{ \rbra{ \sum_{\tilde{j} \in U} D^{i}_{\tilde{i} \tilde{j}} * q_{\tilde{i} \tilde{j}} + (\psi_{UU})_{i \, \tilde{i}} * (1-H_{\tilde{i}}) }^2 } (k) 
\end{split}
\label{} \\
\begin{split}
=& \sum_{\tilde{j} \in U } 
\sbra{ \rbra{ D^{i}_{\tilde{i}\tilde{j}} - (\Psi_{UU})_{i \, \tilde{i}} }^2 * q_{\tilde{i} \tilde{j}} } (k) 
\\
&- 
\sbra{ \rbra{ \sum_{\tilde{j} \in U} D^{i}_{\tilde{i} \tilde{j}} * q_{\tilde{i} \tilde{j}} - (\psi_{UU})_{i \, \tilde{i}} * H_{\tilde{i}} }^2 } (k) . 
\end{split}
\label{}
\end{align}
Hence we obtain $ V^{\vR}_{ee} = v^{\vR}(i,k) $.

\bibliographystyle{plain}

\end{document}